\documentclass[graybox]{svmult}
\usepackage{float}
\usepackage{tabulary}

\usepackage{makeidx}         
\usepackage{graphicx}        
\usepackage{multicol}        
\usepackage[bottom]{footmisc}

\usepackage{amsmath}
\usepackage{newtxtext}       %
\usepackage{newtxmath}       

\usepackage{enumitem}
\usepackage{color}
\newcounter{comments}
\newcommand{\GB}{Gr\"obner basis }
\newcommand{\GBs}{Gr\"obner bases }

\begin{document}
\title*{The Number of Gr\"obner Bases in Finite Fields}
\author{Anyu Zhang, Brandilyn Stigler}
\institute{Anyu Zhang, Southern Methodist University, \email{anyuz@smu.edu}
\and Brandilyn Stigler, Southern Methodist University \email{bstigler@smu.edu}}
%
%
\maketitle

\abstract{
In the field of algebraic systems biology, the number of minimal polynomial models constructed using  discretized data from an underlying system is related to the number of distinct reduced \GBs for the ideal of the data points.  While the theory of \GBs is extensive, what is missing is a closed form for their number for a given ideal. This work contributes connections between the geometry of data points and the number of \GBs associated to small data sets.  Furthermore we improve an existing upper bound for the number of \GBs specialized for data over a finite field.  
%
}

\section{Introduction}
\label{sec:introduction}



Polynomial systems are ubiquitous across the sciences.  While linear approximations are often desired for computational and analytic feasibility, certain problems may not permit such reductions.  In 1965, Bruno Buchberger introduced Gr\"obner bases, which are multivariate nonlinear generalizations of echelon forms \cite{buchberger-thesis,buchberger-translation}. Since this landmark thesis, the adoption of \GBs has expanded into diverse fields, such as geometry \cite{tsai2016}, imagine processing \cite{lin2004},  oil industry \cite{torrente2009},  quantum field theory \cite{maniatis2007}, and systems biology \cite{laubenbacher2004computational}.  

While working with a \GB (GB) of a system of polynomial equations is just as natural as working with a triangularization of a linear system, their complexity can make them cumbersome with which to work: for a general system, the complexity of Buchberger's Algorithm is doubly exponential in the number of variables \cite{buchberger}. The complexity improves in certain settings, such as systems with finitely many real-valued solutions (\cite{BM} is a classic example, whereas \cite{farr} is a more contemporary example), or solutions over finite fields \cite{just}. Indeed much research has been devoted to improving Buchberger's Algorithm and analyzing the complexity and memory usage in more specialized settings (for example, \cite{eder,m4gb}), and even going beyond traditional ways of working with the theory of Gr\"obner bases \cite{larsson};
however most results are for characteristic 0 fields, such $\mathbb R$ or $\mathbb Q$.  


The goal of our work is to consider the \emph{number} of \GBs for a system of polynomial equations over a finite field (which has positive characteristic and consequently all systems have finitely many solutions).  The motivation comes from the work of \cite{laubenbacher2004computational}, in which the authors  presented an algorithm to reverse engineer a model for a biological network from discretized experimental data and made a connection between the number of distinct \emph{reduced} GBs and the number of (possibly) distinct \emph{minimal} polynomial models.
The number of reduced GBs associated to a data set gives a quantitative measure for how ``underdetermined'' the problem of reverse engineering a model for the underlying biological system is.

The Gr\"obner fan geometrically encapsulates all distinct reduced \GBs \cite{mora}. 
In \cite{fukuda} the authors provided an algorithm to compute all reduced GBs.  When their number is too large for enumeration, the method in \cite{dimitrova-volumes} allows one  to sample from the fan.  Finally in \cite{onn}, the authors provide an upper bound for the number of reduced GBs for systems with finitely many solutions; however this bound is much too large for data over a finite field.  To our knowledge, there is no closed form  for the number of reduced GBs, in particular for systems over finite fields with finitely many solutions.

In this paper we make the following contributions with respect to geometric observations of small data sets and computational analyses: 
\begin{enumerate}
    \item geometric characterization of data associated with different numbers of GBs.
    \item formulas for the number of GBs for small data sets over finite fields.
    \item modified upper bound of the number of GBs in the finite field setting.
\end{enumerate}
In Section \ref{sec:background}, we provide the reader with relevant background definitions and results.  In Section \ref{sec:2pts}, we discuss the connection between the number of distinct reduced \GBs for ideals of two points and the geometry of the points; furthermore, we provide a formula to compute the number of GBs associated to 2-point data sets.  We extend the connection to 3 points in Section \ref{sec:3pts}.  Then in Section \ref{sec:bound}, we consider the general setting of any fixed number of points over any finite field and provide an upper bound. We close with a discussion of possible future directions.  We have verified all of the computations referenced in this work, provided illustrative examples throughout the text, and listed data tables in the Appendix.




\section{Background}
\label{sec:background}

Let $K$ be a finite field of characteristic $p>0$.  We will typically consider the finite field  $\mathbb Z_p=\{0,1,\ldots,p-1\}$, that is the field of remainders of integers upon division by $p$ with modulo-$p$ addition and multiplication.  Let $R=K[x_1,\ldots ,x_n]$ be a polynomial ring over $K$.  Finally let $m$ denote the number of points in a subset of $K^n$. Most definitions in this section are taken from \cite{cox}.


A \emph{monomial order} $\prec$ is a total order on the set of all monomials in $R$ that is closed with respect to multiplication and is a well-ordering. 
The \emph{leading term} of a polynomial $g\in R$ is thus the largest monomial for the chosen monomial ordering, denoted as $LT_\prec(g)$. Also we call $LT_\prec(I)=\langle LT_\prec(g):g\in I\rangle$ the \emph{leading term ideal} for an ideal $I$.
\begin{definition}\label{GB}
Let $\prec$ be a monomial order on $R$ and let $I$ be an ideal in $R$. Then $G\subset I$ is a \emph{Gr\"obner basis} for $I$ with respect to $\prec$ if for all $f\in I$ there exists $g\in G$ such that the leading term $LT_\prec(g)$ divides $LT_\prec(f)$. 
\end{definition}

It is well known that \GBs  exist for every $\prec$ and make multivariate polynomial division well defined in that remainders are unique. 
While there are infinitely many orders, there are only finitely many GBs for a given ideal.  As two orders may result in the same \emph{reduced} GB (that is, leading terms have a coefficient of~1 and do not divide other terms in a GB), this results in an equivalence relation where the leading terms of the representative of each equivalence class can be distinguished (underlined) \cite{cox}.  

In this work all GBs are reduced.

\begin{definition}
The monomials which do not lie in $LT_\prec(I)$ are \emph{standard} with respect to $\prec$; the set of standard monomials for an ideal $I$ is denoted by $SM_\prec(I)$. 
\end{definition}

A set of standard monomials  $SM_\prec(I)$ for a given monomial order  forms a basis for $R/I$ as a vector space over $K$.
It is straightforward to check that standard monomials satisfy the following divisibility property: if $x^\alpha \in SM_\prec (I)$ and $x^\beta$ divides $x^\alpha$, then $x^\beta \in SM_\prec (I)$.  This divisibility property on monomials is equivalent to the following geometric condition on lattice points, defined as a \emph{staircase}.

\begin{definition}
\label{def:stircase}
A set $\lambda\subset \mathbb N^n$ is a \emph{staircase} if  for all $u\in \lambda$, $v\leq u$ implies $v\in \lambda$.
\end{definition}

For $S\subseteq K^n$, we call the set $I(S):=\{h\in R\mid h(s)=0 \, \forall s\in S\}$ of polynomials that vanish on $S$ an \textit{ideal of points}. This ideal is computed via standard algebraic geometry techniques as follows:   $I(S)=\bigcap_{i=1}^{m} \langle x_j-s_{ij}\rangle$, that is, $I(S)$ is the intersection of the polynomials that vanish on each point.  Below is the general algorithm to compute polynomial dynamical systems (PDSs) $f=(f_1,\ldots ,f_n):\mathbb Z_p^n\rightarrow \mathbb Z_p^n$ for a given set of data written using the ideal of the input points \cite{laubenbacher2004computational}. 

\textbf{Strategy:} Given input-output data $V=\{(s_1,t_1),\ldots,(s_m,t_m)\}\subset K^n$, find all PDSs that fit $V$ and select a minimal PDS with respect to polynomial division.
\begin{enumerate}
	\item For each $x_j$, compute one interpolating function $f_j\in R$ such that $f_j(s_i)=t_{ij}$. 
	\item Compute the ideal $I=I(\{s_1,\ldots ,s_m\})$ of the input points.
%
\end{enumerate}

Then the \emph{model space} for $V$ is the set $$f+I:=\{(f_1+h_1,\ldots,f_n + h_n):h_i\in I\}$$ of all PDSs which fit the data in $V$ and where $f=(f_1,\ldots,f_n)$ is computed in Step~1. A PDS can be selected from $f+I$ by choosing a monomial order $\prec$, computing a \GB $G$ for $I$, and then computing the remainder (\emph{normal form}) of each $f_i$ by dividing by the polynomials in $G$.  We call 
$$(f_1 \mod G,f_2 \mod G,\ldots,f_n \mod G)$$
the \emph{minimal} PDS with respect to $\prec$, where $G$ is a GB for $I$ with respect to $\prec$. Changing the monomial order may change the resulting minimal PDS.  While it is possible for two reduced GBs to give rise to the same normal form (see \cite{laubenbacher2004computational}), it is still the case that in general a set of data points may have \textit{many} GBs associated to it.  In this way, the number of distinct reduced GBs of $I$ gives an upper bound for the number of different minimal PDSs. Therefore, we aim to find the number of distinct reduced GBs for a given data set. 

\begin{example}
\label{example1}
Consider two inputs $V=\{(0,0),(1,1)\}\in \left(\mathbb Z_2\right)^2$. The corresponding ideal~$I$ of the points in $V$ has 2 distinct reduced Gr\"obner bases, namely 
\begin{center}
$G_1=\{\underline{x_1}-x_2,\underline{x_2^2}-x_2\}, G_2=\{\underline{x_2}-x_1,\underline{x_1^2}-x_1\}$
\end{center}
Here, '\_' marks the leading terms of polynomials in the GBs.
There are 2 resulting minimal models:
any minimal PDS with respect to $G_1$ will be in terms of $x_2$ only as all $x_1$'s are divided out,
while any minimal PDS with respect to $G_2$ will be in terms of $x_1$ only as all $x_2$'s are divided out.
Instead if the inputs are  $\{(0,0),(0,1)\}$, then $I$ has a unique GB 
$\{\underline{x^2_2}-x_2,\underline{x_1}\}$, resulting in a unique minimal PDS.
\end{example}

It is the polynomial $f=x_1-x_2$ that has different leading terms for different monomial orders.  In fact, for monomial orders with $x_1\prec x_2$, the leading term of $f$ will be $x_1$, while for orders with $x_2\prec x_1$ the opposite will be true. We say that $f$ has \textit{ambiguous} leading terms. We will  mark only ambiguous leading terms.

The elements of the quotient ring $R/I$ are equivalence classes of functions defined over the inputs $\{s_1,\ldots s_m\}$.  When working over a finite field, classic results in algebraic geometry state that when the number $m$ of input points is finite, $m=\dim_K R/I$. Since a set of standard monomials is a basis for $R/I$, it follows that each reduced polynomial, $f \mod G$, is written in terms of standard monomials.

Next we state a result about data sets and their complements. 
\begin{theorem}[\cite{robbiano-unique-gb}]
\label{thm:symmetric}
 Let $I$ be the ideal of input points $S$, and let $I^c$ be ideal of  the complement $K^n\setminus S$ of $S$. Then we have $SM_\prec(I) = SM_\prec(I^c)$ and $LT_\prec(I) = LT_\prec(I^c)$ for a given monomial order $\prec$. Hence, we have $\#GB(S) = \#GB(K^n\setminus S)$.
\end{theorem} 

We say that a polynomial $f\in R$ is \emph{factor closed} if every monomial $m\in supp (f)$ is divisible by all monomials smaller than $m$ with respect to an order $\prec$. The following result gives an algebraic description of ideals with unique reduced GBs for any monomial order. 

\begin{theorem}[\cite{robbiano-unique-gb}]%
\label{thm:FC}%
A reduced \GB $G$ with factor-closed generators is reduced for every monomial order; that is, $G$ is the unique reduced \GB for its corresponding ideal.
\end{theorem}

We end this section with a discussion on the number of distinct reduced \GBs for extreme cases.
The set $\mathbb Z_p^n$ contains $p^n$ points. For $n=1$, all ideals have a unique reduced GB since all polynomials are single-variate and as such are factor closed.  We consider cases for $n>1$.
For empty sets  or singletons in $\mathbb Z_p^n$, it is straightforward to show that the ideal of points has a unique reduced GB for any monomial order; that is, for a point $s=(s_1,\ldots ,s_n)$, the associated ideal of $s$ is $I=\langle x_1-s_1,\ldots,x_n-s_n\rangle $ whose generators form a GB and hence is unique (via Theorem \ref{thm:FC}). According to Theorem \ref{thm:symmetric}, the same applies to $p^n-1$ points.  In the rest of this work, we consider the number of reduced GBs for an increasing number of points. 

Note that over a finite field, the relation $x^p-x$ always holds.
\section{Data Sets with $m=2$ Points}
\label{sec:2pts}

Geometric descriptions of data sets can reveal essential features in the underlying network. 
In the problem of counting the number of Gr\"obner bases associated to data sets, the geometric properties of data giving rise to  unique or multiple GBs can provide researchers with a more intuitive way to explore data sets of interest.

\subsection{Two Points over Different Finite Fields}
\label{sec:2ptsgeo}
In this section, we focus on data sets containing two points in $\mathbb Z_p^n$.
\begin{figure}[ht]
\centering
\includegraphics[scale=0.5]{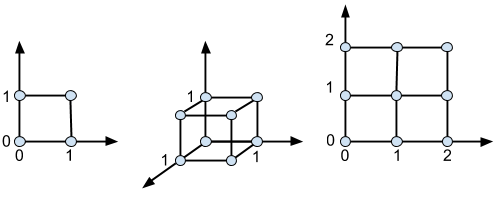}
\caption{The lattice of points in $\mathbb Z_2^2$ (left), in $\mathbb Z_2^3$ (middle), and in $\mathbb Z_3^2$ (right).}
\label{fig:grids} 
\end{figure}

\begin{figure}[ht]
\centering
\includegraphics[scale=0.4]{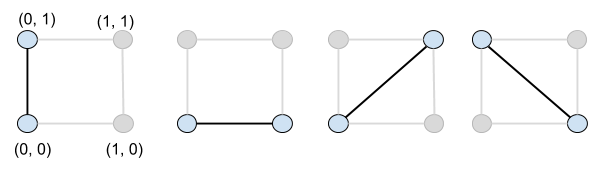}
\caption{Four configurations of pairs of points in $\mathbb Z_2^2$. From left to right: $\{(1,0), (0,1)\}$ and $\{(0,0), (1,0)\}$ each have 1 GB, while $\{(0,0), (1,1)\}$ and $\{(1,0), (0,1)\}$ have 2 distinct GBs.
}
\label{fig:p2n2m2} 
\end{figure}

Consider two coordinates ($n=2$). The left graph in Figure \ref{fig:grids} is the plot of all  points in $\mathbb Z_2^2$.  By decomposing the 2-square on which they lie, we find that pairs of points that lie along horizontal lines will have unique reduced \GBs for any monomial order; see Figure \ref{fig:p2n2m2}.  For example, the set of points $\{(0,0),(0,1)\}$ has ideal of points $\langle  x_1,x_2^2-x_2\rangle$.  Again we can use Theorem \ref{thm:FC} to see that the generators of~$I$ form a unique reduced GB.  Similarly the set of points $\{(1,0),(1,1)\}$ has ideal of points $\langle  x_1-1,x_2^2-x_2\rangle$, which also has a unique reduced GB.  Note that while they have different GBs, they have the same leading term ideal, namely, $\langle  x_1,x_2^2\rangle$. In the same way, pairs of points that lie along vertical lines have unique reduced GBs: sets $\{(0,0),(1,0)\}$ and $\{(0,1),(1,1)\}$ have the unique leading term ideal $\langle  x_1^2,x_2\rangle$.

However pairs of points that lie on diagonals have 2 distinct reduced GBs. For example, the set of points $\{(0,0),(1,1)\}$ has GBs $\{  \underline{x_1}-x_2,x_2^2-x_2\}$ and $\{x_1^2-x_1, \underline{x_2}-x_1\}$ with leading term ideals $\langle x_1, x_2^2\rangle $ and $\langle x_1^2, x_2\rangle $ respectively. 
Similarly the set of points $\{(0,1),(1,0)\}$ has GBs $\{  \underline{x_1}-x_2-1,x_2^2-x_2\}$ and $\{x_1^2-x_1, \underline{x_2}-x_1-1\}$ with leading term ideals $\langle x_1, x_2^2\rangle $ and $\langle x_1^2, x_2\rangle $ respectively. 
%
%
\begin{figure}[H]
\centering
\includegraphics[scale=0.4]{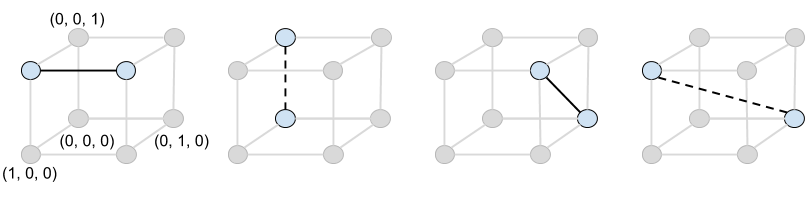}
\caption{Four configurations of pairs of points in $\mathbb Z_2^3$.  From left to right:  
$\{(1,0,1), (1,1,1)\}$ and 
$\{(0,0,0), (0,0,1)\}$ have 1  GB,  
$\{(1,1,1), (0,1,0)\}$ has 2  GBs, and 
$\{(1,0,1), (0,1,0)\}$ has  3 GBs.  
}
\label{fig:p2n3m2} 
\end{figure}
Now consider three coordinates ($n=3$).  The middle graph in Figure \ref{fig:grids} is the plot of all  points in $\mathbb Z_2^3$. In Figure \ref{fig:p2n3m2}, pairs of points that lie on edges of the 3-cube have~1 reduced Gr\"obner basis: for example the set $\{(1,0,1), (1,1,1)\}$ (first from the left in Figure \ref{fig:p2n3m2}) has the unique reduced GB $\{x_1-1, x_2^2-x_2, x_3-1\}$  and  $\{(0,0,0), (0,0,1)\}$ (second) has the associated GB
$\{x_1,x_2,x_3^2-x_3\}$. 
Points that lie on faces of 3-cube have 2 GBs: the set $\{(1,1,1), (0,1,0)\}$ (third) has GBs $\{\underline{x_1}-x_3, x_2-1, x_3^2-x_3\}$ and  $\{x_1^2-x_1, x_2-1, \underline{x_3}-x_1\}$.  Finally points that lie on lines through the interior have 3 GBs:  $\{(1,0,1), (0,1,0)\}$ (fourth) has GBs $\{\underline{x_1}-x_3, \underline{x_2}-x_3-1, x_3^2-x_3\}$, $\{\underline{x_1}-x_2-1, x_2^2-x_2, \underline{x_3}-x_2-1\}$, and $\{x_1^2-x_1, \underline{x_2}-x_1-1, \underline{x_3}+x_1\}$. From the summary of 2 points in Figure \ref{fig:p2n2m2}, we know that a unique GB arises when the points lie on horizontal or vertical edges.  Furthermore the simultaneous change of  both variables (coordinates)  will lead to nonunique GBs. This behavior is reasonable, as in this case, the data cannot distinguish which variable is the leading variable. 
As $n$ increases, we can see the general trend; that is, the number of distinct GBs coincides with the number of coordinate changes between the two points.


\begin{figure}[ht]
\centering
\includegraphics[scale=0.3]{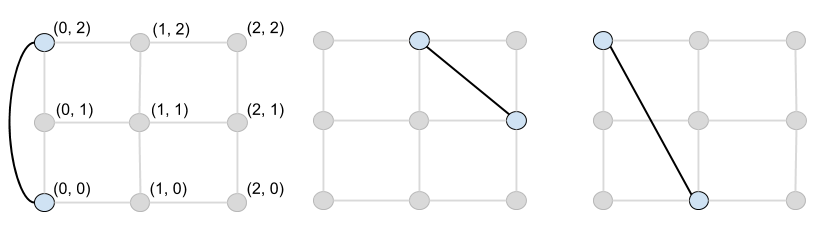}
\caption{Three configurations of points in $\mathbb Z_3^2$. From left to right: $\{(0,0), (0,2)\}$ has 1 GB, while 
$\{(1,2), (2,1)\}$ and
$\{(0,2), (1,0)\}$ each have 2 distinct GBs.  
}
\label{fig:p3n2m2} 
\end{figure}

Lastly, consider non-Boolean fields. Let $p=3$ and $n=2$. The right graph in Figure~\ref{fig:grids} is the plot of all  points in $\mathbb Z_3^2$.  
Similar to the Boolean case, pairs of points that lie on horizontal or vertical lines have one associated reduced \GB for any monomial order, while pairs of points that lie on any skew line have two distinct GBs; see Figure \ref{fig:p3n2m2}.  
For example, the set of points $\{(0,0),(0,2)\}$ has ideal of points $\langle  x_1,x_2^2+x_2\rangle$, which has a unique reduced GB via Theorem \ref{thm:FC}.  On the other hand, the set of points $\{(1,2),(2,1)\}$ has two GBs, namely 
$\{  \underline{x_1}+x_2,x_2^2+1\}$ and 
$\{  x_1^2-1, \underline{x_2}+x_1,\}$ with leading term ideals
$\langle x_1, x_2^2\rangle $ and $\langle x_1^2, x_2\rangle $ respectively.  

As the number $n$ of coordinates increases, the number of coordinate changes between the two points determines the number of distinct reduced GBs; however increasing the number $p$ of states does not affect the number of reduced GBs, as we will see in the next section.

\subsection{General Formula for Two Points}
We conclude with a theorem which summarizes our findings for sets of two points with any number of variables and any number of states. Let $m$ be the number of points in $K^n$. Define $N_n^m$ to be the number of \GBs  for ideals with $m$ points in $K^n$.

\begin{definition} 
\label{def:2pointspiece}
For a given coordinate value $x$, define the following piece-wise function. 

 \begin{eqnarray}
   B_0(x) = \left\{
     \begin{array}{lr}
       1 & : x =0\\
       0 & : x >0
     \end{array}
   \right.
\end{eqnarray}
\end{definition}

\begin{theorem}
\label{thm:2points}
Let $P=(p_1,\ldots,p_n), Q=(q_1,\ldots,q_n)\in K^n$ and let $I$ be the ideal of the points $P,Q$.  The number of distinct reduced \GBs for $I$  is given by
\begin{eqnarray}
\label{formula:2pts}
N^2_n =n-\sum_{i=1}^n B_0\left(|p_i-q_i|\right).
\end{eqnarray}
\end{theorem}

\smallskip
\begin{proof}
Let $P=(p_1, \ldots, p_n)$ and $Q=(q_1,\ldots, q_n)$ be two points in $K^n$.  
Recall that the ideal $I$ of the points can be computed via intersections:
\begin{center}
$I=\langle x_1-p_1,\ldots,x_n-p_n\rangle \cap \langle x_1-q_1,\ldots,x_n-q_n\rangle$.
\end{center}
%
Organizing the multiplication of terms in a matrix, we get the following:
\[A=
\begin{bmatrix} 
(x_1 -p_1)(x_1-q_1) & (x_1-p_1)(x_2-q_2) & \ldots &(x_1-p_1)(x_n-q_n)\\
(x_2-p_2)(x_1-q_1) & (x_2 -p_2)(x_2-q_2) & \ldots &(x_2-p_2)(x_n-q_n)\\
 \vdots& \vdots & \ddots &\vdots\\
 (x_n-p_n)(x_1-q_1) & (x_n-p_n)(x_2-q_2) & \ldots & (x_n -p_n)(x_n-q_n)
\end{bmatrix}
\]
%
In general, subtracting polynomials $A_{ij} = (x_i-q_i)(x_j-p_j)$  with the same leading term $A_{ji}=(x_i-q_j)(x_j-p_i)$ yields  $$A_{ij}-A_{ji}=(q_j-p_j)x_i+(p_i-q_i)x_j+q_ip_j-p_iq_j$$ which has a smaller-power leading term $x_i$ or $x_j$. 
%
The system $A^*$ of simplified equations after subtracting the polynomials below the diagonal from the ones above the diagonal is
\begin{align}
\label{d1}
x^2_1-(p_1+q_1)x_1+p_1q_1,\, A_{12}-A_{21},\,
 A_{13}-A_{31},
\nonumber \ldots
, A_{1n}-A_{n1} \tag{$A^*_1$}\\
\label{d2}
x^2_2-(p_2+q_2)x_2+p_2q_2, \,
\nonumber A_{23}-A_{32},\ldots, 
\nonumber A_{2n}-A_{n2}, \tag{$A^*_2$} \\ 
\nonumber \vdots\\ 
\label{dn}
x^2_n-(p_n+q_n)x_n+p_nq_n \tag{$A^*_n$}
\end{align}
Notice that the expressions $$x^2_1-(p_1+q_1)x_1+p_1q_1, \, x^2_2-(p_2+q_2)x_2+p_2q_2, \dots, x^2_n-(p_n+q_n)x_n+p_nq_n$$ in the rows (\ref{d1}), (\ref{d2}), and (\ref{dn}), corresponding to the diagonal entries in the matrix $A$, are univariate and have the same leading term for any monomial order.

Suppose $x_1$ is the smallest variable in some order. Then for all forms $x_1+x_i$, the leading term is $x_i$, for $i=2,3,\ldots , n$. Hence, the set of leading terms is  $\{x^2_1,x_2,x_3,\ldots,x_n\}$. With $n$ free choices of smallest variable $x_k$ to fix the leading term of $A_{ij}-A_{ij}$, for $k=1, \ldots, n$, we can get at most $n$ different GBs.
However, the coordinates of $P$ and $Q$ will affect the leading terms. If one coordinate is the same, say $p_1 = q_1$,
then the expressions in \ref{d1} are of the form $A_{1j} - A_{j1} = (p_j-q_j)x_1 + p_1q_j-q_1p_j$ for $j=1,\ldots ,n$. The leading term of these polynomials is $x_1$. Then there are $n-1$ free choices left for smallest variable $x_k$ for $k = 2, \ldots, n$. In this case there will be at most $n-1$ GBs. By iteration, Formula \ref{formula:2pts} can be proved. The number of distinct reduced GBs can be decreased with more coordinates being equal. At last, if $n-1$ coordinates are equal, two points generate a unique GB for any monomial order.
\end{proof}

\begin{example}
Let $P=(1,0,0), Q=(0, 1, 0)$, based on Formula \ref{formula:2pts}, we can get number of Gr\"obner bases of the ideal of these two points:
$$N^3_2 = 3 - (B_0(|p_1-q_1|) + B_0(|p_2-q_2|) + B_0(|p_3-q_3|)= 3 -(0+0+1) = 2.$$
The distinct GBs for the ideal of the points are
$$\{x_3, x_2^2+x_2, \underline{x_1}+x_2+1\},
\{x_3, x_1+\underline{x_2}+1, x_1^2+x_1\}.$$

\end{example}

\begin{corollary}
The maximum number of distinct reduced \GBs for an ideal of two points in $\mathbb Z_p^n$ is $n$.
\end{corollary}

These results hold for data sets of two points over any finite field $\mathbb Z_p$ in any number of coordinates $n$. 
%
\section{Data Sets with $m \geq 3$ Points}
\label{sec:3pts}

Here we extend the ideas of Section \ref{sec:2pts} to sets with more points and offer geometric observations for a small number of choices for $p$ states and $n$ coordinates.

\subsection{Three Points over Different Finite Fields}
\label{sec:p2m3}

We begin by considering the Boolean base field. 

Let  $n=2$. For three points in $\mathbb Z_2^2$ all ideals have a unique reduced \GB via Theorem \ref{thm:symmetric} and the fact that ideals of a single point have only one reduced GB for any monomial order.

Now let $n=3$. 

\begin{figure}[ht]
\centering
\includegraphics[scale=.5]{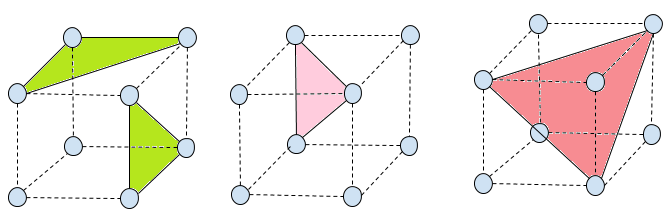}
\caption{Configurations of sets of 3 points in $\mathbb Z_2^3$ corresponding to different numbers of GBs.  Points that are in configurations similar to the green triangles (left) have a unique reduced \GB for any monomial order; the pink triangle (middle) has two distinct GBs; and the red triangle (right) has three distinct GBs.}
\label{fig:p2n3}       
\end{figure}

\begin{example}
Consider the point configurations in Figure \ref{fig:p2n3}.  The data set corresponding to the  green triangle on the top ``lid'' of the leftmost 3-cube is $S_1 = \{(0, 1, 0),  (0, 1, 1), (1, 1, 1)\}$ and has  a unique associated reduced Gr\"obner basis:
$$\{x_3^2+x_3, x_2+1, x_1x_3+x_1, x_1^2+x_1\}.$$
The data set corresponding to the pink triangle in middle 3-cube is $S_2 = \{(0, 0, 1), (0, 1, 1), (1, 1, 0)\}$ and has   two associated reduced GBs:
$$ \{x_3^2+x_3, x_2x_3+x_2+x_3+1, x_2^2+x_2, \underline{x_1}+x_3+1\}; \{x_1+\underline{x_3}+1, x_2^2+x_2, x_1x_2+x_1, x_1^2+x_1\}.$$
Finally the data set corresponding to the  red triangle in the rightmost  3-cube is $S_3 = \{(1, 0, 0), (0, 1, 0), (1, 1, 1)\}$ and has  three associated reduced GBs:
 $$ \{x_3^2+x_3, x_2x_3+x_3, x_2^2+x_2, \underline{x_1}+x_2+x_3+1\},$$ $$\{x_3^2+x_3, x_1+\underline{x_2}+x_3+1, x_1x_3+x_3, x_1^2+x_1\},$$ $$ \{x_1+x_2+\underline{x_3}+1, x_2^2+x_2, x_1x_2+x_1+x_2+1, x_1^2+x_1\}. $$
\end{example}

The example illustrates that points that lie on faces of the 3-cube have 1 Gr\"obner basis; points forming a triangle which lies in the interior with 2 collinear vertices have 2 distinct GBs, and points in other configurations have 3 GBs. Based on the characteristics in Boolean fields , we provide two formulas for the number of \GB $N_2^3$ and $N_3^3$ for sets of three points in $\mathbb Z_2^2$ and $\mathbb Z_2^3$. Recall that~$N_n^m$ is the number of GBs for ideals with $m$ points in $K^n$. Motivated by Theorem \ref{thm:2points}, we construct the following piece-wise functions. 

\begin{definition} 
Let $x\in \mathbb Z$.  Define the following functions: 
\label{def:3pointspiece}
\begin{align}
   B_1(x) = \left\{
     \begin{array}{lr}
       1 & : x =1\\
       0 & : x =2
     \end{array}
   \right.
\qquad\qquad & 
   B_2(x) = \left\{
     \begin{array}{lr}
       0 & : x <0\\
       x & : x \geq{0}
     \end{array}
   \right.
\end{align}
\end{definition}

Let $ P=(p_1,p_2),Q=(q_1,q_2),R=(r_1,r_2)\in \mathbb Z_2^2$. Set
\begin{center}
$(s_1,s_2,s_3)=(N^2_2(P,Q),N^2_2(P,R),N^2_2(Q,R))$.
\end{center}
The number of distinct reduced \GBs for ideals of 3 points in $\mathbb Z_2^2$ is 
\begin{eqnarray}
N^3_2 =n-B_2\left(\sum_{i=1}^3 B_1\left(s_i\right)-1\right). 
\label{formula:3pts}
\end{eqnarray}

Let $ P=(p_1,p_2, p_3),Q=(q_1,q_2, q_3),R=(r_1,r_2, r_3)\in \mathbb Z_2^3$. Set
\begin{center}
$(w_1,w_2,w_3)=(N^2_3(P,Q),N^2_3(P,R),N^2_3(Q,R)).$
\end{center}
The number of distinct reduced \GBs for ideals of 3 points in $\mathbb Z_2^3$ is 
\begin{eqnarray}
N^3_3 =n-\sum_{i=1}^3 
B_1\left(w_i\right).
\label{formula1:3pts}
\end{eqnarray}
\begin{example} 
Let $P=(1,0), Q=(0, 1), R = (0, 0)\in \mathbb Z_2^2$. Using Formula \ref{formula:3pts}, the number of distinct reduced Gr\"obner bases of the ideal of these three points:$$N^3_2 = 2 - B_2(B_1(N_2^2(P, Q)) + B_1(N_2^2(P, R)) + B_1(N_2^2(Q, R)) - 1)$$ $$= 2 - B_2(B_1(2)+B_1(1)+B_1(1)-1) = 2-B_2(0+1+1-1)=2-1=1.$$
In fact, $\{\{x_1^2+x_1, x_1x_2, x_2^2+x_2\}\}$ is the unique reduced GB for the ideal of points.
\end{example}

\begin{example}
Let $P=(1,0,0), Q=(0, 1, 0), R=(0, 0, 1)\in \mathbb Z_2^3$.  Using Formula \ref{formula1:3pts}, we compute the number of Gr\"obner bases of the ideal of these three points:$$N^3_3 = 2 - (B_1(N_3^2(P, Q)) + B_1(N_3^2(P, R)) + B_1(N_3^2(Q, R)))$$ $$ = 3 -(B_1(3-1) + B_1(3-1) + B_1(3-1))= 3 - (0+0+0) = 3.$$
The distinct GBs for the ideal of the points are 
$$\{x_3^2+x_3, x_2x_3, x_2^2+x_2, \underline{x_1}+x_2+x_3+1\}$$
$$\{x_3^2+x_3, x_1+\underline{x_2}+x_3+1, x_1x_3, x_1^2+x_1\} $$
$$\{x_1+x_2+\underline{x_3}+1, x_2^2+x_2, x_1x_2, x_1^2+x_1\}.$$
\end{example}

Formulas \ref{formula:3pts} and \ref{formula1:3pts} for three points in $\mathbb Z_2^n$ for $n=2,3$ are motivated by the calculation of the number of distinct  \GBs for two points  in $\mathbb Z_2^n$ for any $n$; see Equation \ref{formula:2pts}.
A more general formula with larger $n$ or larger $p$ (non-Boolean fields) is hard to generate. Especially, the construction of piece-wise functions becomes harder and unpredictable, with the increasing number of variables and number of states. 
%
%
%


\medskip

Now we turn our attention to non-Boolean base fields.  Let $p=3$ and $n=2$.  See the right graph in Figure \ref{fig:grids} for a plot of all points in $\mathbb Z_3^2$. 
\begin{figure}[H]
\centering
\includegraphics[scale=0.3]{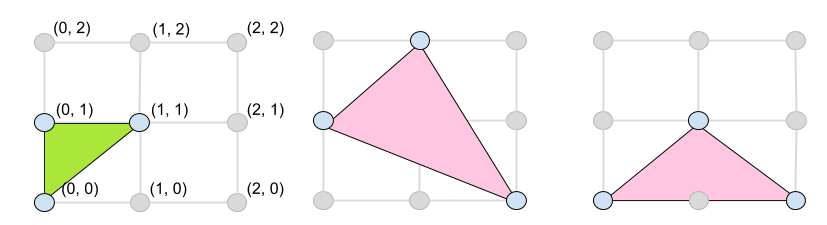}
\caption{
Configurations of sets of 3 points in $\mathbb Z_3^2$ corresponding to unique and non-unique Gr\"obner bases.  Points that are in configurations similar to the green triangles (left) have a unique reduced \GB for any monomial order; the pink triangles (middle and right) have two distinct GBs.}
\label{fig:p3n2m3}       
\end{figure}
\begin{example}
Consider the point configurations in Figure \ref{fig:p3n2m3}. The data set corresponding to the green triangle (left) is $S_1 = \{(0, 0), (0, 1),  (1, 1)\}$ and has a unique associated reduced Gr\"obner basis:
$$\{x_2^2-x_2, x_1x_2-x_1, x_1^2-x_1\}.$$
The data set corresponding to the pink triangle (middle) is
$S_2={(0, 1),(1, 2),(2, 0)}$ and has  two associated reduced GBs:
$$\{x_2^3-x_2, \underline{x_1}-x_2+1\}, \{-x_1+\underline{x_2}-1, x_1^3-x_1\}.$$
The data set corresponding to the pink triangle (right) is
$S_3={(0, 1),(1, 2),(2, 0)}$ and has two associated reduced GBs:
$$\{x_2^3-x_2, x_1x_2^2-x_1x_2+x_2^2-x_2, x_1^2-x_1x_2+x_1-x_2\}, \{x_2^3-x_2, -x_1^2+x_1x_2-x_1+x_2, x_1^3-x_1\}.$$
\end{example}

In Figure \ref{fig:p3n2m3}, we see that 3 points that lie on a line or form a green triangle with vertices distance 1 from each other have unique Gr\"obner bases, while 3 points that form other configurations in red triangles have 2 distinct GBs. Similar to the results with 2 points in Figure \ref{fig:p2n2m2}, if only one variable (coordinate) changes (configurations of vertical or horizontal lines), unique \GB be will generated. In Figure~\ref{fig:p2n3} in $\mathbb Z_2^3$ (Boolean case) and Figure~\ref{fig:p3n2m3} in $\mathbb Z_3^2$ (non-Boolean case), with an increase in the number of diagonal edges, data sets in these configurations will generate non-unique Gr\"obner bases.

\begin{figure}[H]
\includegraphics[scale=.55]{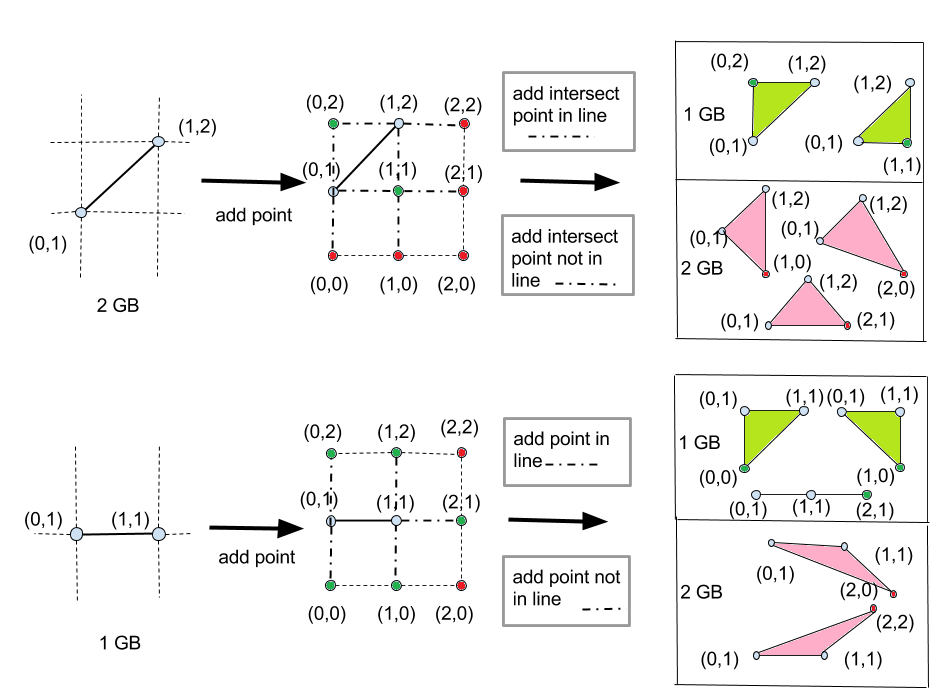}
\caption{Adding a point to a two-point set.  In the top panel, the initial blue points do not have a unique GB. Adding a green point reduces the number of GBs to unique one, while adding points in red positions will not reduce the number of GBs. In the bottom panel, the initial points have a unique GB. Adding a green point keeps the unique number of GBs, while adding a red point increases the number of GBs. }

\label{fig:addpts} 
\end{figure}

To generalize the geometric pattern from small data sets to larger data sets, we start with configurations of 2 points then add a point to characterize the observed numbers of \GBs for 3 points. Using Figure \ref{fig:addpts}, adding a green point on horizontal or vertical lines will decrease the number of GBs, while adding a red point on diagonal lines will not result in a unique GB. 
%
%

Based on the geometric characteristics associated with unique GBs, in the next section we state a conjecture for decreasing the number of \GB by adding points in so-called \textit{linked} positions.

\subsection{Larger Numbers of Points over Different Finite Fields}
\begin{definition}
\label{def:linkedposition}
Given a set $S$ of points, we say that a point $q$ is in a \textit{linked} position with respect to the points in $S$ if $q$ lies on the same grid lines as the points in $S$.  
\end{definition}
For example, the green points in Figure \ref{fig:addpts} are in linked position with respect to the blue points.

\begin{conjecture}
\label{conjecture:addpoints} 
Let $S$ be a set of points, $q$ is a point not in $S$, and $T = S\cup \{q\}$. 
If $q$ is in a linked position and the convex hull of the points in $T$ does not contain ``holes'' (\textit{i.e}, lattice points not in $T$), then $\#GB(T) \leq \#GB(S)$. 
\end{conjecture}
For $n=2,3$ variables and $p= 2,3$ states, we constructed all possible subsets of points in $\mathbb Z_p^n$ and computed the number of GBs for data sets up to 6 points.  While these are modest results, they are relevant in the sense that experimental data tend to be small (less than 10 input-output observations).
\begin{figure}[H]
\centering\includegraphics[scale=.6]{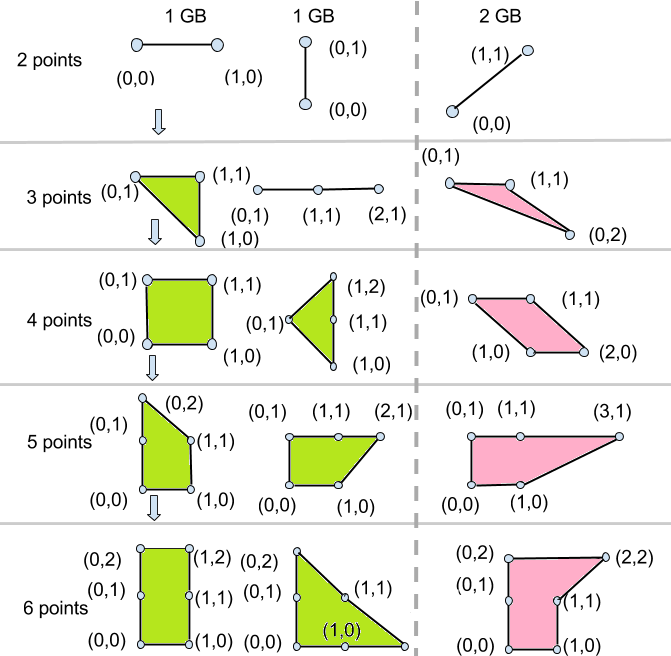}
\caption{Point configurations based on the number of GBs for $m=2,\ldots ,6$. The left two columns contain points that form green triangles and correspond to a unique Gr\"obner basis. The right  column contains the pink triangles corresponding to non-unique GBs.}
\label{fig:morepts}       
\end{figure}

In Figure \ref{fig:morepts}, the points in the configuration of green triangles that have a uniquely associated GB.
By adding more points in linked positions, the augmented data set will keep the unique \GB like first column (moving downward) in Figure \ref{fig:morepts}.
Based on the geometric characteristics in the above pictures, we can summarize the following rules to aid researchers decrease the number of candidate models as enumerated by the number of distinct GBs:
\begin{enumerate}
\item  \textit{For two points},  fewer changing coordinates in the data points will lead to fewer GBs.  In the simplest case, if only one coordinate changes, a unique model will be generated.
\item \textit{For three points}, more points lying on horizontal or vertical edges will reduce the number of GBs.  A unique GB arises when the data lie on a horizontal line, a vertical line or form a right triangle. 
\item \textit{In the process of adding points}, if researchers want to decrease or keep the number of minimal models, the better candidates of new data points are those in linked positions with respect to an existing data set: this guarantees more points lying on horizontal or vertical edges.
\end{enumerate}

By adding points in linked positions, data sets with multiple \GBs can be transformed to data sets with unique GB, as the following example suggests. 

\begin{example}
Consider data sets in $\mathbb Z_2^4$.
Let $S_{max}$ be a data set whose ideal of points has the maximum number of GBs. Define  $S_{unique}=S_{max}\cup S_{add}$ where $S_{add}$ is a collection of points such that the augmented data set $S_{unique}$ has an ideal of points with a unique GB. The table summarizes for different sized sets how many points must be added to an existing data set to guarantee a unique GB.

\begin{center}
\begin{tabular}{|c|ccccccccccccc|}
\hline
 $\max(\#GBs)$ & 4 & 5 & 6 & 13 & 12 & 13 & 9 & 13 & 12 & 13 & 6 & 5 & 4 \\ \hline
 $|S_{max}|$ & 2 & 3 & 4 & 5 & 6 & 7 & 8 & 9 & 10 & 11 & 12 & 13 & 14 \\ \hline
 $|S_{unique}|$ & 5 & 5 & 8 & 11 & 11 & 11 & 11 & 12 & 15 & 15 & 15 & 15 & 15  \\ \hline
 $|S_{add}|$ & 3 & 2 & 4 & 6 & 5 & 4 & 3 & 3 & 5 & 4 & 3 & 2 & 1 \\  \hline
\end{tabular}
\end{center}

\end{example}

\section{Upper Bound for the Number of Gr\"obner Bases}
\label{sec:bound}

We now focus on the general setting of subsets of any size $m$ in $\mathbb Z_p^n$ for any $p$ and~$n$.  In \cite{onn} the authors proved that an upper bound for the maximum number of GBs for an ideal of $m$ points in $F^n$ over an arbitrary base field $F$ is 
\begin{equation}
    \max \#GBs(m,n) = m^{2n{\frac{n-1}{n+1}}}
    \label{onnbound}
\end{equation}
%
where
\begin{enumerate}
\item the number of vertices of any lattice polytope $\mathcal P$ is $O\left(vol(\mathcal P)^\frac{n-1}{n+1}\right)$ \cite{andrews}
\item distinct reduced GBs can be identified with the vertices of a specially constructed polytope (generalized polygon) of volume $(m^2)^n$; see right panel in Figure 5.
\end{enumerate}
For finite fields with $p$ states, this bound becomes unnecessarily large for even small~$m$. Since polytopes in a finite field are contained in a hypercube of volume $p^n$, we aim to modify the above result accordingly.
\begin{figure}
\begin{center}
\includegraphics[width=8cm]{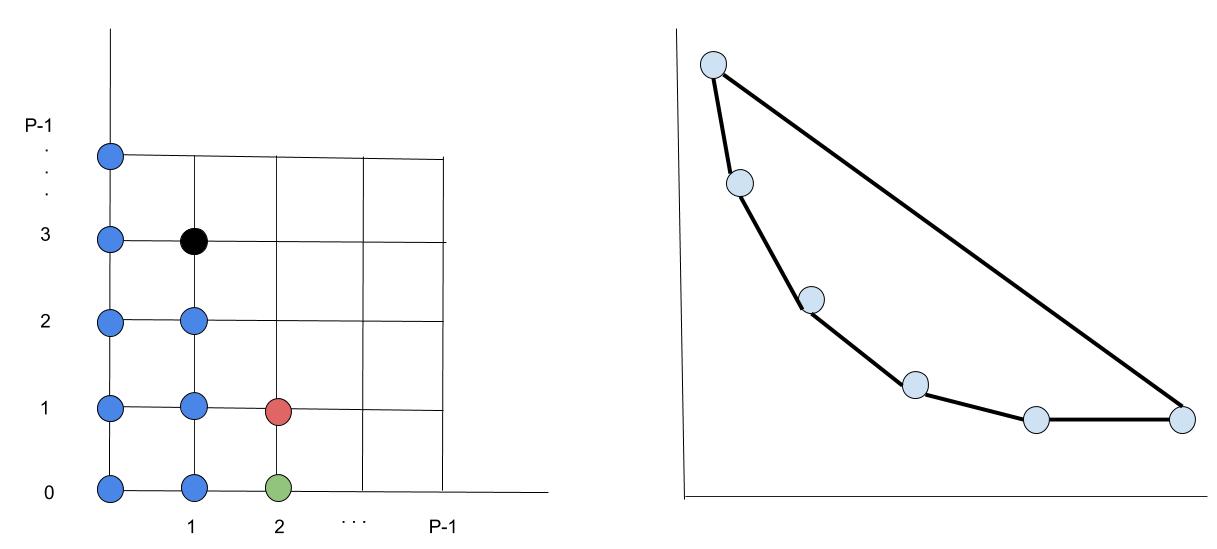}
\caption{The lattice graph of a staircase and corresponding polytope as defined in \cite{onn}.}
\end{center}
\label{fig:polyandstaircase}
\end{figure}
As the authors in  \cite{onn} showed that the sum of the coordinates of a staircase of $m$ points (see left panel in Figure \ref{fig:polyandstaircase}) corresponds to a vertex of a certain polytope, we must therefore count the number of ways to place $m$ points on the lattice.

Suppose $r$ blue points have been placed (see Figure \ref{fig:polyandstaircase}).  We wish to count the number of ways to place to next point.  The red point violates the staircase property.  The only choice is the green or black point.  Note that the black point maximizes the sum of the coordinates.

The polytope $\mathcal P$ is contained in the convex body $n$-simplex 
\begin{eqnarray}
conv\{0,\max \sum_{i=1}^{m}\lambda_1 e_1,\max\sum_{i=1}^{m}\lambda_2 e_2,\ldots,\max\sum_{i=1}^{m}\lambda_n e_n\}
\end{eqnarray}
For an infinite number of states, $\max \sum_{i=1}^{m}\lambda_{i}=\binom m2$.  However,
for $p<\infty$ states,
\begin{eqnarray}
\max \sum_{i=1}^{m}\lambda_{i}={\frac{p(p-1)}{2}}\left \lfloor m/p \right \rfloor +\dbinom {m\mod p}{2}
\end{eqnarray}
Modifying the bound in \cite{onn} for $p$ states gives
\begin{align}
vol(\mathcal P) &\leq \frac{1}{n!}\left({\max \sum_{i=1}^{m}\lambda }\right)^n \\
 &\leq \frac{1}{n!}\left({\frac{p(p-1)}{2}\left \lfloor m/p \right \rfloor +\frac{(m\mod p)(m\mod p-1)}{2} }\right)^n \\ 
 &= O\left(\left(p^2\left \lfloor m/p \right \rfloor +({m\mod p})^2\right)^n\right).
  \label{equ:upperbound1}
\end{align}

When a data set is empty, plugging in $m=0$ into Equation \ref{equ:upperbound1} results in a calculated upper bound of 0; similarly for the case of choosing $p^n$ points. (Note that the bound in \cite{onn} makes the same calculation.)  So, the equation applies for nonempty sets of size less than $p^n$. For completion sake, we include the extreme cases $m = 0 $ and $m=p^n$ into Equation  \ref{equ:upperbound1}.

\begin{theorem}
\label{thm:upperbound}
The number of distinct reduced \GBs for an ideal of $m$ points in $\mathbb Z_p^n$ is at most
\begin{align}
\label{equ:upperbound2}
   U(n, m, p) = \left\{
     \begin{array}{lr}
       \left(p^2 \left \lfloor m/p \right \rfloor +(m\mod p)^2\right)^{n\frac{n-1}{n+1}} & : \lfloor p^n/2 \rfloor\geq m > 0\\
       1 & : m = 0
     \end{array}
   \right.
\end{align} 
When $m > \lfloor p^n/2\rfloor$, then the number of \GBs is given by $U(n, p^n-m, p)$. 
\end{theorem}

Based on Equation \ref{equ:upperbound2}, we apply the result in Theorem \ref{thm:symmetric} and summarize the modified formula for number of Gr\"obner basis as follows:
\begin{align}
 \label{equ:upperbound4}
   U(n, m, p) = \left\{
     \begin{array}{lr}
     \left(p^2 \left \lfloor m/p \right \rfloor +(m\mod p)^2\right)^{n\frac{n-1}{n+1}} & : \lfloor p^n/2 \rfloor\geq m > 0\\
       \left(p^2 \left \lfloor (p^n-m)/p \right \rfloor +((p^n-m)\mod p)^2\right)^{n\frac{n-1}{n+1}} & : p^n > m > \lfloor p^n/2 \rfloor\\
       1 & : m = 0, p^n
     \end{array}
   \right.
\end{align}
It is straightforward to show that our bound grows much slower than the bound  $O\left(m^{2n\frac{n-1}{n+1}}\right)$ reported in \cite{onn}, which we have also verified computationally. Below is a table of selected numerical results of the new upper bound in comparison to the values of the original upper bound in \cite{onn}.

\begin{figure}[H]\centering
   \begin{minipage}{0.48\textwidth}
     \frame{\includegraphics[width=1\linewidth]{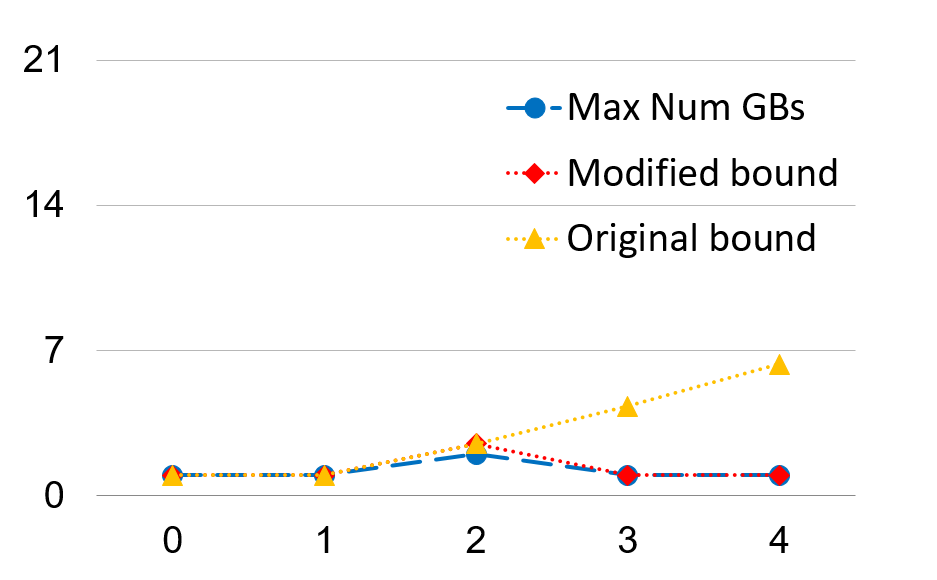}}
     \caption{$p=2$, $n=2$.}\label{fig:p2n2bound}
   \end{minipage}
   \begin{minipage}{0.48\textwidth}
     \frame{\includegraphics[width=1\linewidth]{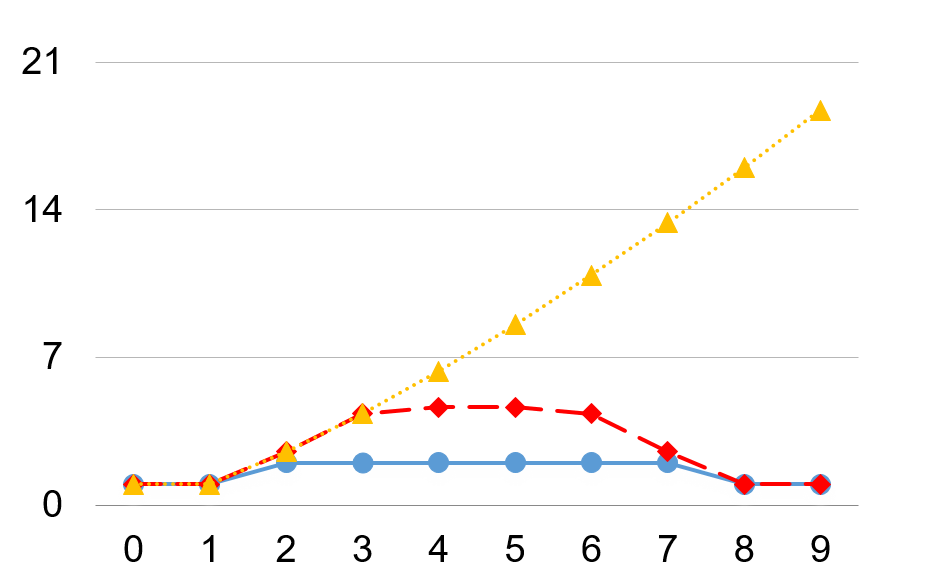}}
     \caption{$p=3$, $n=2$.}\label{fig:p3n2bound}
   \end{minipage}
    \begin{minipage}{0.48\textwidth}
    \vspace{0.4cm}
     \frame{\includegraphics[width=1\linewidth]{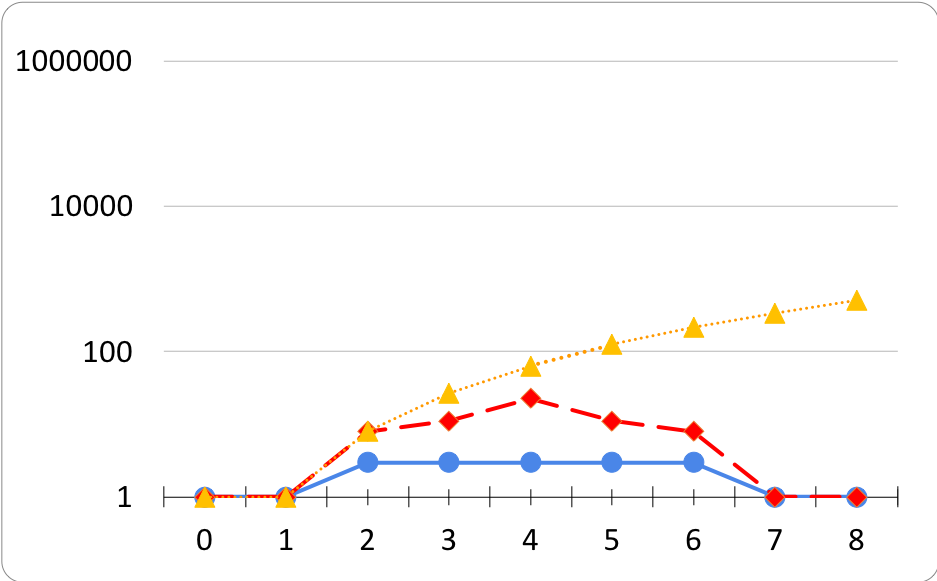}}
     \caption{$p=2$, $n=3$.}\label{fig:p2n3bound}
   \end{minipage}
   \begin{minipage}{0.48\textwidth}
   \vspace{0.4cm}
     \frame{\includegraphics[width=1\linewidth]{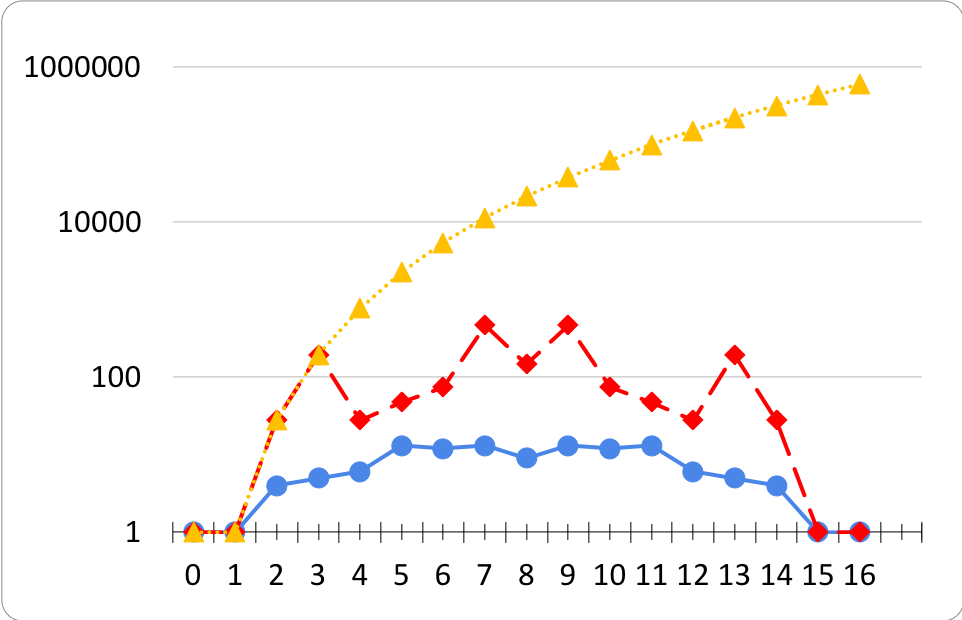}}
     \caption{$p=2$, $n=4$.}\label{fig:p2n4bound}
   \end{minipage}
   \caption{Plots comparing the maximum number of Gr\"obner bases. The caption in each plot indicates the values of $p$ and $n$ for $\mathbb Z_p^n$.  In each case, all subsets of size $m$ are computed, where $m$ ranges from 0 to $p^n$ and listed on the horizontal axis.  The vertical axis is the maximum number of GBs for a set with of size $m$. The blue solid line with dots shows the actual maximum number of GBs.  The yellow dotted line with triangles is the number predicted by Equation \ref{onnbound}, where the red dashed line with squares is the number predicted by Equation \ref{equ:upperbound4}.}
\end{figure}

The experiments of the original bound and the modified bound are compared in the following case: $n=4$ (four variables) and $p=2$ (Boolean states). We compared the formula results with the maximum number of GBs in each case.
With the changing of the number of states from $p=2$ (Figure \ref{fig:p2n2bound}) to $p=3$ (Figure \ref{fig:p3n2bound}), the number of GBs will increase, and the difference between original bound and actual bound becomes larger. Our modified bound's performance is much better than the original bound, especially in cases with a large number of points.




Considering the effects of the increasing number of variables $n$, the number of variables $n$ decides the power value in original bound in Equation \ref{onnbound} and in the modified bound in Equation \ref{equ:upperbound4}. Hence, changing the number of states only from 2 to 3 will lead to large differences in predictions. For example, for the case $n=4$, $p=2$ and number of points $m = 5$ in Figure \ref{fig:p2n4bound}, the original bound is over 2000, while the modified bound is much closer to the maximum number of GBs. The original bound  will not be helpful for researchers to estimate in finite field settings because the difference increases rapidly away from the actual number  with an increasing number of points.  However, the modified bound  provides reasonable estimates of maximum number of GBs in different finite fields shown from Figure \ref{fig:p2n2bound} to Figure \ref{fig:p2n4bound}.
%

\section{Discussion}
\label{sec:discussion}
This work relates the geometric configuration of data points with the number of associated Gr\"obner bases.  In particular we provided some insights into which configurations lead to unique GBs. We give formulas for the specific number of \GBs for small data sets, and also provide researchers with a way to decrease the number of GBs by adding extra points in so-called \textit{linked} positions. At last, we developed a modified upper bound specialized for finite fields, which has been tested in in a variety of cases.  An implication of this work is a more computationally accurate way to predict the number of distinct minimal models which may aid researchers in estimating the computational cost before running physical experiments. 
 %
%

Increasing $p$, $n$ or $m$ will all inflate the difference between the predicted number of GBs and the actual number. The performance of the modified bound works well with large $p$ and $m$. However, based on Table \ref{tab:upperboundcompare4} in the Appendix, the modified bound, though better than original bound, still has large difference from actual values for $n>4$.
Hence, how to decrease the effects of the increasing number of variables is future work for upper bound estimation.



\section*{Appendix}
\label{sec:appendix}
Below we provide tables summarizing the comparison of the maximum number of distinct reduced \GBs to the predictions made by the original bound listed in Equation \ref{onnbound} and the modified bound listed in Equation \ref{equ:upperbound4}.
The second column shows the actual maximum number as computed for all sets in $\mathbb Z_p^n$ of size given in the first column.  The third column gives the values computing using Equation \ref{onnbound} and the last column gives the values computing using Equation \ref{equ:upperbound4}.
%
\begin{table}[H]
\caption{$p=2, n=2$}
\label{tab:upperboundcompare0}
\begin{center}
\begin{tabular}{|c|c|c|c|}
\hline
\# of points & max \# of GBs & original bound & modified bound  \\
\hline
0 & 1 & 1 & 1 \\
1 & 1 & 1 & 1  \\
2 & 2 & 3 & 3 \\
3 & 1 & 4 & 1 \\
4 & 1 & 6 & 1 \\
\hline
\end{tabular}
\end{center}
\end{table}
\vspace{-0.5in}
\begin{table}[H]
\caption{$p=2, n=3$}
\label{tab:upperboundcompare1}
\begin{center}
\begin{tabular}{|c|c|c|c|}
\hline
\# of points & max \# of GBs & original bound & modified bound  \\
\hline
0 & 1 & 1 & 1\\
1 & 1 & 1 & 1  \\
2 & 3 & 8 & 8 \\
3 & 3 & 27 & 11 \\
4 & 3 & 64 & 23 \\
5 & 3 & 125 & 11 \\
6 & 3 & 216 & 8 \\
7 & 1 & 343 & 1  \\
8 & 1 & 512 & 1 \\
\hline
\end{tabular}
\end{center}
\end{table}
\vspace{-0.5in}
\begin{table}[H]
\caption{$p=2, n=4$}
\label{tab:upperboundcompare2}
\begin{center}
\begin{tabular}{|c|c|c|c|}
\hline
\# of points & max \# of GBs & original bound & modified bound  \\
\hline
1 & 1 & 1 & 1  \\
2 & 4 & 28 & 28  \\
3 & 5 & 195 & 195 \\
4 & 6 & 776 & 28 \\
5 & 13 & 2.26E+03 & 48 \\
6 & 12 & 5.43E+03 & 74 \\
7 & 13 & 1.14E+04 & 471  \\
8 & 9 & 1.93E+04 & 147 \\
\hline
\end{tabular}
\end{center}
\end{table}
\vspace{-0.5in}
\begin{table}[H]
\caption{$p=3, n=2$}
\label{tab:upperboundcompare3}
\begin{center}
\begin{tabular}{|c|c|c|c|}
\hline
\# of points & max \# of GBs & original bound & modified bound  \\
\hline
1 & 1 & 1 & 1  \\
2 & 2 & 3 & 3 \\
3 & 2 & 4 & 4 \\
4 & 2 & 6 & 5 \\
5 & 2 & 9 & 5 \\
6 & 2 & 11 & 4 \\
7 & 2 & 13 & 3  \\
8 & 1 & 16 & 1 \\
9 & 1 & 19 & 1 \\
\hline
\end{tabular}
\end{center}
\end{table}
\vspace{-0.5in}
\begin{table}[H]
\caption{$m = 4$ points and $p = 2$}
\label{tab:upperboundcompare4}
\begin{center}
\begin{tabular}{|c|c|c|c|}
\hline
\# of variables & max \# of GBs & original bound & modified bound  \\
\hline
2 & 1 & 6 &  1 \\
3 & 3 & 64 & 27 \\
4 & 5 & 776 & 147 \\
5 & 8 & 10321& 1024\\
\hline
\end{tabular}
\end{center}
\end{table}
\end{document}